\documentclass[12pt,reqno]{amsart}
\usepackage{cite}
\usepackage{amsfonts}
\usepackage{amssymb, amsmath, latexsym}
\usepackage{longtable}
\usepackage{lipsum}
\usepackage{hyperref}
\usepackage{mleftright}
\usepackage{graphicx,mathtools}
\usepackage{mathtools}

\usepackage{array}

\numberwithin{equation}{section}
\theoremstyle{plain}
\newtheorem{theorem}{Theorem}[section]
\newtheorem{definition}{Definition}[section]

\newtheorem{corollary}[theorem]{Corollary}

\title[ Divisibility of an analogue of $t$-core partition function by powers of primes] {Divisibility of an analogue of $t$-core partition function by powers of primes}

\author[P. Talukdar]{Pranjal Talukdar}
\address{Department of Mathematical Sciences, Tezpur University,  Assam 784028, India}
\email{pranjaltalukdar113@gmail.com}

\allowdisplaybreaks
\begin{document}
\begin{abstract}
A partition of a positive integer $n$ is said to be $t$-core if none of its hook lengths are divisible by $t$.  Recently, two analogues, $\overline{a}_t(n)$ and $\overline{b}_t(n)$, of the $t$-core partition function, $c_t(n)$, have been introduced by Gireesh, Ray and Shivashankar \cite{grs} and Bandyopadhyay and Baruah \cite{bb}, respectively. In this article, we prove the lacunarity of $\overline{b}_t(n)$ modulo arbitrary powers of 2 and 3 for $t=3^\alpha m$ where $\gcd(m,6)$=1. For a fixed positive integer $k$ and prime numbers $p_i\geq 5$, we also study the arithmetic density of $\overline{b}_t(n)$ modulo $p_i^k$ where $t=p_1^{a_1}\cdots p_m^{a_m}$. We further prove an infinite family of congruences for $\overline{b}_3(n)$ modulo arbitrary powers of 2 by employing a result of Ono and Taguchi on the nilpotency of Hecke operators.
\end{abstract}
\maketitle

\vskip 3mm
\noindent {\footnotesize 2010 Mathematical Reviews Classification Numbers:  11P83, 05A17, 11F11.}

\section{Introduction and statement of the results}\label{Introduction}

For complex numbers $a$ and $q$ such that $\mid q\mid <1$, we define the infinite $q$-product as
\begin{align*}
\left(a;q\right)_\infty:=\prod_{j=0}^\infty\left(1-aq^j\right).
\end{align*}
In the sequel, for brevity, we set $f_n:=\left(q^n;q^n\right)_\infty$ for integers $n\ge1$.

An integer partition $\pi=\{\pi_1,\pi_2,\cdots,\pi_k\}$ of a positive integer $n$ is a non-increasing sequence of natural numbers such that $\displaystyle\sum_{i=1}^k \pi_i=n$. We denote the number of partitions of $n$ by $p(n)$. The Ferrers–Young diagram of $\pi$ is an array of nodes with $\pi_i$ nodes in the $i$th row. The $(i,j)$ hook is the set of nodes directly to the right of $(i,j)$ node, together with the set of nodes directly below it, as well as the $(i,j)$ node itself. The hook number, $H(i,j)$, represents the total number of nodes on the $(i,j)$ hook. For $t\geq 2$, a partition of $n$ is called $t$-core if none of the hook numbers are divisible by $t$. An illustration of the Ferrers-Young diagram of the partition $4+3+1$ of 8 with hook numbers is given below:
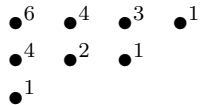
\begin{center}
\begin{tabular}{ cccc }  
 $\bullet^6$ &  $\bullet ^4$&  $\bullet^3$&  $\bullet^1$\\ 
 $ \bullet^4$& $ \bullet^2$ &  $\bullet^1$&\\ 
 $\bullet^1$ &  & & 
\end{tabular}
\end{center}
Clearly, for $t\geq7$, the partition $4+3+1$ of 8 is a $t$-core partition. The number of $t$-core partitions is denoted by $c_t(n)$. In an existence result, Granville and Ono \cite{go} proved that for $t\geq 4$, every natural number $n$ has a $t$-core. Core partitions play an important role in many branches of mathematics. Apart from partition theoretic aspects, they have connections with other areas such as representation theory of symmetric groups and symmetric functions (see \cite{go}, \cite{st}). For a brief survey of the topic, one can refer to \cite{ck}.

Let $G(q):=\displaystyle\sum_{n=0}^{\infty}a(n)q^n$ be an integral power series and $0\leq r < M$. The arithmetic density $\delta_r(G,M;X)$ is defined as
\begin{align*}
\delta_r(G,M;X):=\dfrac{\#\left\{0\leq n\leq X: a(n)\equiv r\pmod M\right\}}{X}.
\end{align*}
An integral power series $G$ is called \textit{lacunary modulo} $M$ if
\begin{align*}
\lim_{X\rightarrow\infty}\delta_0(G,M;X)=1,
\end{align*}
that is, almost all of the coefficients of $G$ are divisible by $M$.

Arithmetic densities of $c_t(n)$ modulo arbitrary powers of 2, 3 and primes greater than or equal to 5 are recently studied by Jindal and Meher \cite{jm}. 

Recall that for $\vert ab\vert <1$, Ramanujan's general theta function $f(a,b)$ is given by
\begin{align*}
f(a,b)=\displaystyle\sum_{n=-\infty}^\infty a^{n(n+1)/2}b^{n(n-1)/2}=(-a;ab)_\infty(-b;ab)_\infty(ab;ab)_\infty,
\end{align*}
where the last equality is Jacobi’s famous triple product identity \cite[p. 35, Entry 19]{bcb3}.

Three special cases of $f(a,b)$ are:
\begin{align}
\varphi(-q)&:=f(-q,-q)=\displaystyle\sum_{n=-\infty}^\infty (-1)^nq^{n^2}=\dfrac{f_1^2}{f_2}\label{p-q},\\
\psi(-q)&:=f(-q,-q^3)=\displaystyle\sum_{n=-\infty}^\infty (-q)^{n(n+1)/2}=\dfrac{f_1f_4}{f_2}\label{s-q},\\
f(-q)&:=f(-q,-q^2)=\displaystyle\sum_{n=-\infty}^\infty (-1)^nq^{n(3n-1)/2}=f_1\label{f-q}.
\end{align}

In the notation of Ramanujan's theta functions, the generating function of $c_t(n)$ is given by
\begin{align}
\displaystyle\sum_{n=0}^{\infty}c_t(n)q^n=\dfrac{f^t(-q^t)}{f(-q)}\label{ct1}.
\end{align}

Recently, Gireesh, Ray and Shivashankar \cite{grs} introduced a new function, $\overline{a}_t(n)$, by substituting $\varphi(-q)$ in place of $f(-q)$ in the generating function of $c_t(n)$ (in \eqref{ct1}), namely
\begin{align}
\displaystyle\sum_{n=0}^{\infty}\overline{a}_t(n)q^n=\dfrac{\varphi^t(-q^t)}{\varphi(-q)}=\dfrac{f_2f_t^{2t}}{f_1^2f_{2t}^t}\label{a t}.
\end{align}
Using Ramanujan's theta functions and $q$-series, they deduced certain multiplicative formulas and arithmetic identities for $\overline{a}_t(n)$ for $t=$ 2, 3, 4 and 8. They also studied the divisibility of $\overline{a}_t(n)$ modulo arbitrary powers of primes greater than 5 using the theory of modular forms. Their result can be stated as in the following theorem.
\begin{theorem}
Let $t=p_1^{a_1}\cdots p_m^{a_m}$ where $p_i$’s are prime numbers greater than or equal to 5. Then for every positive integer $j$, the set
\begin{align*}
\left\{ n\in \mathbb{N}: \overline{a}_t(n)\equiv 0\pmod{p_i^j}\right\}
\end{align*}
has arithmetic density 1.
\end{theorem}
Additionally, they proved a Ramanujan type congruence for $\overline{a}_5(n)$ modulo 5 by using an algorithm developed by Radu and Sellers \cite{rs}. Bandyopadhyay and Baruah \cite{bb} established new identities connecting $\overline{a}_5(n)$ and $c_5(n)$. They also deduced a reccurence relation for $\overline{a}_5(n)$. The author \cite{pt} studied the arithmetic density of $\overline{a}_{t}(n)$ modulo arbitrary powers of 2 and 3 where $t=3^\alpha m$. The density results in \cite{pt} can be stated as in the following theorem.

\begin{theorem}\label{a 2k}
Let $k\geq 1$, $\alpha\geq 0$ and $m \geq 1$ be integers with \textup{gcd}$(m, 6)=1$. For $p\in\{2,3\}$, the
set
\begin{align*}
\left\{n \in \mathbb{N} : \overline{a}_{3^\alpha m}(n) \equiv 0 \pmod {p^k}\right\}
\end{align*}
has arithmetic density 1.
\end{theorem}

Very recently, Bandyopadhyay and Baruah \cite{bb} considered another analogue $\overline{b}_t(n)$ of $t$-core, which is defined by
\begin{align}
\displaystyle\sum_{n=0}^{\infty}\overline{b}_t(n)q^n=\dfrac{\psi^t(-q^t)}{\psi(-q)}=\dfrac{f_2f_t^tf_{4t}^t}{f_1f_4f_{2t}^t}\label{b t}.
\end{align}
In particular, they studied the function $\overline{b}_5(n)$ and deduced some new identities connecting $c_5(n)$, $\overline{a}_5(n)$ and $\overline{b}_5(n)$. They also proved some reccurence relations and vanishing coefficients results for $\overline{b}_5(n)$.

A result containing certain sufficient conditions for lacunarity of products of eta-quotients modulo arbitrary powers of primes were recently proved by Cotron \textit{et. al.} \cite[Theorem 1.1]{cot}. We observe that the eta-quotients associated with $\overline{b}_t(n)$ do not satisfy these conditions, which makes the study of arithmetic density of $\overline{b}_t(n)$ more interesting. In this article, we prove the lacunarity of $\overline{b}_t(n)$ modulo arbitrary powers of 2 and 3 when $t=3^\alpha m$. To be precise, we prove the following results.

\begin{theorem}\label{b 2k}
Let $k\geq 1$ be a fixed positive integer. Then for $\alpha\geq 0$ and $m \geq 1$ with \textup{gcd}$(m, 6)=1$, $\overline{b}_{3^\alpha m}(n)$ is almost always divisible by $2^k$, namely,
\begin{align*}
\lim_{X\rightarrow\infty}\dfrac{\#\{0\leq n\leq X:\overline{b}_{3^\alpha m}(n)\equiv 0\pmod{2^k}\}•}{X}=1.
\end{align*}
\end{theorem}

\begin{theorem}\label{b 3k}
Let $k\geq 1$ be a fixed positive integer. Then for $\alpha\geq 0$ and $m \geq 1$ with \textup{gcd}$(m, 6)=1$, $\overline{b}_{3^\alpha m}(n)$ is almost always divisible by $3^k$, namely,
\begin{align*}
\lim_{X\rightarrow\infty}\dfrac{\#\{0\leq n\leq X:\overline{b}_{3^\alpha m}(n)\equiv 0\pmod{3^k}\}•}{X}=1.
\end{align*}
\end{theorem}

We also study the density of $\overline{b}_t(n)$ modulo arbitrary powers of primes greater than or equal to 5 for certain general values of $t$. In fact, we prove the following result.
\begin{theorem}\label{b pk}
Let $k\geq 1$ be a fixed positive integer and $a_1,a_2, \cdots, a_m$ be non negative integers for some positive integer $m$, then for $t = p_1^{a_1}\cdots p_m^{a_m}$ where $p_i$’s are prime numbers greater
than or equal to 5, $\overline{b}_{t}(n)$ is almost always divisible by $p_i^k$, namely,  
\begin{align*}
\lim_{X\rightarrow\infty}\dfrac{\#\{0\leq n\leq X:\overline{b}_t(n)\equiv 0\pmod{p_i^k}\}}{X}=1.
\end{align*}
\end{theorem}

As a consequence of the above theorem, we obtain the following divisibility result for $\overline{b}_{p}(n)$.
\begin{corollary}
Let $k\geq 1$ be a fixed positive integer and $p\geq 5$ is a prime number. Then $\overline{b}_{p}(n)$ is almost always divisible by $p^k$, namely,  
\begin{align*}
\lim_{X\rightarrow\infty}\dfrac{\#\{0\leq n\leq X:\overline{b}_p(n)\equiv 0\pmod{p^k}\}}{X}=1.
\end{align*}
\end{corollary}

The fact that the action of Hecke algebras on spaces of modular forms of level 1 modulo 2 is locally nilpotent was first observed by Serre and proved by Tate (see \cite{ser1}, \cite{ser2}, \cite{ta}). Their result was later generalized to higher levels by Ono and Taguchi \cite{ot}. This result has been used for proving congruences for certain partition functions. (For example see \cite{a}, \cite{bss}, \cite{jm1}, \cite{sb} etc.) In this paper, we observe that the
eta-quotient associated to $\overline{b}_3(n)$ is a modular form whose level is in
the list of Ono and Taguchi. Thus, we use a result of Ono and Taguchi to prove the following congruences for $\overline{b}_3(n)$.

\begin{theorem}\label{cong2}
Let $n$ be a nonnegative integer. Then there exists an integer $u\geq 0$
such that for every $v\geq 1$ and distinct primes $q_1,\ldots,q_{u+v}$ coprime to 6, we have
\begin{align*}
\overline{b}_3\left(\dfrac{q_1\cdots q_{u+v}\cdot n-24}{24}\right)\equiv 0\pmod{2^v}
\end{align*}
whenever $n$ is coprime to $q_1,\ldots,q_{u+v}$.
\end{theorem}

We organize the paper in the following way. In Section \ref{prelim}, we state some preliminaries of the theory of modular forms and Hecke operators. Then we prove Theorems \ref{b 2k}-\ref{cong2} using the properties of modular forms in Sections \ref{pr b2k}-\ref{pr cong2}. We end the paper with some concluding remarks in Section \ref{cr}.
\section{Preliminaries}\label{prelim}

In this section, we recall some basic facts and definitions on modular forms. For details, we refer the readers to \cite{kob} and \cite{ono}.

We first define the matrix groups
\begin{align*}
\textup{SL}_2(\mathbb{Z})&:=\left\{\begin{bmatrix}
a & b \\
c & d 
\end{bmatrix} : a, b, c, d \in \mathbb{Z}, ad-bc=1\right\},\\
\Gamma_0(N)&:=\left\{\begin{bmatrix}
a & b \\
c & d 
\end{bmatrix} \in \textup{SL}_2(\mathbb{Z}) : c\equiv 0 \pmod N\right\},\\
\Gamma_1(N)&:=\left\{\begin{bmatrix}
a & b \\
c & d 
\end{bmatrix} \in 
\Gamma_0(N) : a\equiv d\equiv 1 \pmod N\right\},
\end{align*}
and
\begin{align*}
\Gamma(N):=\left\{\begin{bmatrix}
a & b \\
c & d 
\end{bmatrix} \in 
\textup{SL}_2(\mathbb{Z}) : a\equiv d\equiv 1 \pmod N, \textup{ and } b\equiv c\equiv 0 \pmod N \right\},
\end{align*}
where $N$ is a positive integer. A subgroup $\Gamma$ of $\textup{SL}_2(\mathbb{Z})$ is called a congruence subgroup if $\Gamma(N)\subseteq \Gamma$ for some $N$ and the smallest $N$ with this property is called the level of $\Gamma$. For instance, $\Gamma_0(N)$ and $\Gamma_1(N)$ are congruence subgroups of level $N$.

Let $\mathbb{H}$ denote the upper half of the complex plane. The group 
\begin{align*}
\textup{GL}_2^{+}(\mathbb{R}):=\left\{\begin{bmatrix}
a & b \\
c & d 
\end{bmatrix} \in 
\textup{SL}_2(\mathbb{Z}) : a, b, c, d \in \mathbb{R} \textup{ and } ad-bc>0 \right\}
\end{align*}
acts on $\mathbb{H}$ by $\begin{bmatrix}
a & b \\
c & d 
\end{bmatrix} z = \dfrac{az+b}{cz+d}$. We identify $\infty$ with $\dfrac{1}{0}$ and define $\begin{bmatrix}
a & b \\
c & d 
\end{bmatrix} \dfrac{r}{s} = \dfrac{ar+bs}{cr+ds}$, where $\dfrac{r}{s} \in \mathbb{Q} \cup \{ \infty \}$. This gives an action of $\textup{GL}_2^{+}(\mathbb{R})$ on the extended upper half plane $\mathbb{H}^*=\mathbb{H}\cup \mathbb{Q}\cup \{\infty\}$. Suppose that $\Gamma$ is a congruence subgroup of $\textup{SL}_2(\mathbb{Z})$. A cusp of $\Gamma$ is an equivalence class in $\mathbb{P}^1=\mathbb{Q}\cup\{\infty\}$ under the action of $\Gamma$. 

The group $\textup{GL}_2^{+}(\mathbb{R})$ also acts on functions $f:\mathbb{H}\rightarrow \mathbb{C}$. In particular, suppose that $\gamma=\begin{bmatrix}
a & b \\
c & d 
\end{bmatrix}$ $\in \textup{GL}_2^{+}(\mathbb{R})$. If $f(z)$ is a meromorphic function on $\mathbb{H}$ and $\ell$ is an integer, then define the slash operator $|_{\ell}$ by
\begin{align*}
(f|_{\ell}\gamma)(z):=(\textup{det}(\gamma))^{\ell/2}(cz+d)^{-\ell}f(\gamma z).
\end{align*}

\begin{definition}
Let $\Gamma$ be a congruence subgroup of level $N$. A holomorphic function $f:\mathbb{H}\rightarrow \mathbb{C}$ is called a modular form with integer weight $\ell$ on $\Gamma$ if the following hold:
\begin{enumerate}
\item[(1)] We have
\begin{align*}
f\left(\dfrac{az+b}{cz+d}\right)=(cz+d)^{\ell}f(z)
\end{align*}
for all $z \in \mathbb{H}$ and and all $\begin{bmatrix}
a & b \\
c & d 
\end{bmatrix}$ $\in \Gamma$.
\item[(2)] If $\gamma \in$ $\textup{SL}_2(\mathbb{Z})$, then $(f|_{\ell}\gamma)(z)$ has a Fourier expansion of the form 
\begin{align*}
(f|_{\ell}\gamma)(z)=\displaystyle\sum_{n\geq 0}a_{\gamma}(n)q_{N}^n,
\end{align*}
where $q:=e^{2\pi iz/N}$.
\end{enumerate}
\end{definition}
For a positive integer $\ell$, the complex vector space of modular forms of weight $\ell$ with respect to a congruence subgroup $\Gamma$ is denoted by $M_{\ell}(\Gamma)$.
\begin{definition}\cite[Definition 1.15]{ono}
If $\chi$ is a Dirichlet character modulo $N$, then we say that a modular form $f \in M_{\ell}(\Gamma_{1}(N))$ has Nebentypus character $\chi$ if 
\begin{align*}
f\left(\dfrac{az+b}{cz+d}\right)=\chi(d)(cz+d)^{\ell}f(z)
\end{align*}
for all $z\in \mathbb{H}$ and all $\begin{bmatrix}
a & b \\
c & d 
\end{bmatrix}$ $\in \Gamma_{0}(N)$. The space of such modular forms is denoted by $M_{\ell}(\Gamma_{0}(N),\chi)$.
\end{definition}

The relevant modular forms for the results of this paper arise from eta-quotients. Recall that the Dedekind eta-function $\eta(z)$ is defined by
\begin{align*}
\eta(z):=q^{1/24}(q;q)_{\infty}=q^{1/24}\displaystyle\prod_{n=1}^{\infty}(1-q^n),
\end{align*}
where $q:=e^{2\pi iz}$ and $z\in \mathbb{H}$. A function $f(z)$ is called an eta-quotient if it is of the form 
\begin{align*}
f(z)=\displaystyle\prod_{\delta|N}\eta(\delta z)^{r_{\delta}},
\end{align*}
where $N$ is a positive integer and $r_{\delta}$ is an integer. Now, we recall two important theorems from \cite[p. 18]{ono} which
will be used later.

\begin{theorem}\cite[Theorem 1.64]{ono}\label{level}
If $f(z)=\displaystyle\prod_{\delta|N}\eta(\delta z)^{r_{\delta}}$ is an eta-quotient such that $\ell=\dfrac{1}{2}\displaystyle\sum_{\delta|N}r_{\delta}$ $\in \mathbb{Z}$,
\begin{align*}
\displaystyle\sum_{\delta|N}\delta r_{\delta}\equiv 0 \pmod {24} \quad\quad \text{and}\quad\quad
\displaystyle\sum_{\delta|N}\dfrac{N}{\delta} r_{\delta}\equiv 0 \pmod {24},
\end{align*}
then $f(z)$ satisfies 
\begin{align*}
f\left(\dfrac{az+b}{cz+d}\right)=\chi(d)(cz+d)^{\ell}f(z)
\end{align*}
for every $\begin{bmatrix}
a & b \\
c & d 
\end{bmatrix}$ $\in \Gamma_{0}(N)$. Here the character $\chi$ is defined by $\chi(d)$ $:=$ $\left(\dfrac{(-1)^{\ell}s}{d}\right),$ where $s:=$ $\displaystyle\prod_{\delta|N}\delta^{r_{\delta}}$.
\end{theorem}

Consider $f$ to be an eta-quotient which satisfies the conditions of Theorem \ref{level} and that the associated weight $\ell$ is a positive integer. If $f(z)$ is holomorphic at all the cusps of $\Gamma_{0}(N)$, then $f(z) \in M_{\ell}\left(\Gamma_{0}(N), \chi\right)$. The necessary criterion for determining orders of an
eta-quotient at cusps is given by the following theorem.
\begin{theorem}\cite[Theorem 1.64]{ono}\label{order}
Let $c$, $d$ and $N$ be positive integers with $d|N$ and gcd$(c,d)$=1. If $f$ is an eta-quotient satisfying the conditions of Theorem \ref{level} for $N$, then the order of vanishing of $f(z)$ at the cusp $(c/d)$ is 
\begin{align*}
\dfrac{N}{24}\displaystyle\sum_{\delta|N}\dfrac{\text{gcd}(d,\delta)^2r_{\delta}}{\text{gcd}(d,N/d)d\delta}.
\end{align*}
\end{theorem}

We now recall a deep theorem of Serre \cite[Page 43]{ono} which will be used in proving Theorems \ref{b 2k}--\ref{b pk}.
\begin{theorem} \cite[p. 43]{ono}\label{Ser}
Let $g(z) \in M_\ell(\Gamma_0(N),\chi)$ has Fourier expansion
\begin{align*}
g(z)=\displaystyle\sum_{n=0}^\infty b(n)q^n \in \mathbb{Z}[[q]].
\end{align*}
Then for a positive integer $r$, there is a constant $\alpha > 0$ such that
\begin{align*}
\#\{0<n\leq X: b(n)\not\equiv 0 \pmod r\}=\mathcal{O}\left(\dfrac{X}{(log X)^\alpha}\right).
\end{align*}
Equivalently
\begin{align*}
\displaystyle{\lim_{X \to \infty}} \dfrac{\#\{0<n\leq X: b(n)\not\equiv 0 \pmod r\}}{X}=0.
\end{align*}
\end{theorem}

Finally, we recall the definition of Hecke operators.  Let $m$ be a positive integer and $f(z)=\displaystyle\sum_{n=0}^{\infty}a(n)q^n$ $\in M_{\ell}(\Gamma_{0}(N),\chi)$. Then the action of Hecke operator $T_{m}$ on $f(z)$ is defined by
\begin{align*}
f(z)|T_m:=\displaystyle\sum_{n=0}^{\infty}\left(\displaystyle\sum_{d|\text{gcd}(n,m)}\chi(d)d^{\ell-1}a\left(\dfrac{nm}{d^2}\right)\right)q^n.
\end{align*}
In particular, if $m=p$ is prime, then we have
\begin{align}
\label{operator}f(z)|T_p:=\displaystyle\sum_{n=0}^{\infty}\left(a(pn)+\chi(p)p^{\ell-1}a\left(\dfrac{n}{p}\right)\right)q^n.
\end{align}
We note that $a(n)=0$ unless $n$ is a nonnegative integer.
\section{Proof of Theorem \ref{b 2k}}\label{pr b2k}
Putting $t=3^{\alpha}m$ in \eqref{b t}, we have
\begin{align}
\displaystyle\sum_{n=0}^{\infty}\overline{b}_{3^{\alpha}m}(n)q^n=\dfrac{f_2f_{3^{\alpha}m}^{3^{\alpha}m}f_{4\cdot3^{\alpha}m}^{3^{\alpha}m}}{f_1f_4f_{2\cdot3^{\alpha}m}^{3^{\alpha}m}}.\label{b 3am}
\end{align}
We define
\begin{align*}
E_{\alpha,m}(z):=\dfrac{\eta^2\left(2^43^{\alpha+1}mz\right)}{\eta\left(2^53^{\alpha+1}mz\right)}.
\end{align*}

Applying the binomial theorem, for any integer $k\geq 1$, we have
\begin{align}
E^{2^k}_{\alpha,m}(z)=\dfrac{\eta^{2^{k+1}}\left(2^43^{\alpha+1}mz\right)}{\eta^{2^k}\left(2^53^{\alpha+1}mz\right)}\equiv 1\pmod{2^{k+1}}\label{E 2k}.
\end{align}

Next, we define 
\begin{align*}
F_{\alpha,m,k}(z)&:=\dfrac{\eta(48z)\eta^{3^{\alpha}m}\left(2^33^{\alpha+1}mz\right)\eta^{3^{\alpha}m}\left(2^53^{\alpha+1}mz\right)}{\eta(24z)\eta(96z)\eta^{3^{\alpha}m}\left(2^43^{\alpha+1}mz\right)}E^{2^k}_{\alpha,m}(z)\notag\\
&=\dfrac{\eta(48z)\eta^{3^{\alpha}m}\left(2^33^{\alpha+1}mz\right)\eta^{2^{k+1}-3^{\alpha}m}\left(2^43^{\alpha+1}mz\right)}{\eta(24z)\eta(96z)\eta^{2^k-3^{\alpha}m}\left(2^53^{\alpha+1}mz\right)}.
\end{align*}

Using \eqref{b 3am} and \eqref{E 2k}, we obtain
\begin{align}
F_{\alpha,m,k}(z)&\equiv\dfrac{\eta(48z)\eta^{3^{\alpha}m}\left(2^33^{\alpha+1}mz\right)\eta^{3^{\alpha}m}\left(2^53^{\alpha+1}mz\right)}{\eta(24z)\eta(96z)\eta^{3^{\alpha}m}\left(2^43^{\alpha+1}mz\right)}\notag\\
&\equiv q^{3\left(3^{2\alpha}m^2-1\right)}\dfrac{f_{48}f_{2^33^{\alpha+1}m}^{3^{\alpha}m}f_{2^53^{\alpha+1}m}^{3^{\alpha}m}}{f_{24}f_{96}f_{2^43^{\alpha+1}m}^{3^{\alpha}m}}\notag\\
&\equiv\displaystyle\sum_{n=0}^{\infty}\overline{b}_{3^{\alpha}m}(n)q^{24n+3\left(3^{2\alpha}m^2-1\right)}\pmod{2^{k+1}}\label{F amk}.
\end{align}

Now, we will show that $F_{\alpha,m,k}(z)$ is a modular form. Applying Theorem \ref{level}, we find that the level of $F_{\alpha,m,k}(z)$ is $N=2^53^{\alpha+1}mM$, where $M$ is the smallest positive integer such that
\begin{align*}
2^53^{\alpha+1}mM\left(\dfrac{-1}{24}+\dfrac{1}{48}+\dfrac{-1}{96}+\dfrac{3^{\alpha}m}{2^33^{\alpha+1}m}+\dfrac{2^{k+1}-3^{\alpha}m}{2^43^{\alpha+1}m}+\dfrac{-2^{k}+3^{\alpha}m}{2^53^{\alpha+1}m}\right)\equiv 0\pmod{24},
\end{align*}
which implies 
\begin{align*}
3\cdot2^kM\equiv 0\pmod{24}.
\end{align*}
Therefore, $M=4$ and the level of $F_{\alpha,m,k}(z)$ is $N=2^7 3^{\alpha+1}m$.

The cusps of $\Gamma_0\left(2^73^{\alpha+1}m\right)$ are given by fractions $c/d$ where $d|2^73^{\alpha+1}m$ and $\gcd(c, d) = 1$ (see for example \cite[page 5]{ono1}). By Theorem \ref{order}, $F_{\alpha,m,k}(z)$ is holomorphic at a cusp $c/d$ if and only if
\begin{align*}
&- \dfrac{\gcd (d,24)^2}{24}+\dfrac{\gcd (d,48)^2}{48}-\dfrac{\gcd (d,96)^2}{96}+3^\alpha m \dfrac{\gcd \left(d,2^33^{\alpha+1} m\right)^2}{2^33^{\alpha+1} m}\\
&+\left(2^{k+1}- 3^\alpha m\right) \dfrac{\gcd \left(d,2^43^{\alpha+1} m\right)^2}{2^43^{\alpha+1} m}-\left(2^k-3^\alpha m\right) \dfrac{\gcd \left(d,2^53^{\alpha+1} m\right)^2}{2^53^{\alpha+1} m}\geq 0.
\end{align*}
Equivalently, $F_{\alpha,m,k}(z)$ is holomorphic at a cusp $c/d$ if and only if
\begin{align*}
L:=3^{\alpha}m(-4G_1+2G_2-G_3+4G_4-2G_5+1)+2^k(4G_5-1)\geq 0,
\end{align*}
where $G_1=\dfrac{\gcd (d,24)^2}{\gcd \left(d,2^5 3^{\alpha+1} m\right)^2}$, $G_2=\dfrac{\gcd (d,48)^2}{\gcd \left(d,2^5 3^{\alpha+1} m\right)^2}$, $G_3=\dfrac{\gcd (d,96)^2}{\gcd \left(d,2^5 3^{\alpha+1} m\right)^2}$,

 $G_4=\dfrac{\gcd (d,2^3 3^{\alpha+1} m)^2}{\gcd \left(d,2^5 3^{\alpha+1} m\right)^2}$ and $G_5=\dfrac{\gcd (d,2^4 3^{\alpha+1} m)^2}{\gcd \left(d,2^5 3^{\alpha+1} m\right)^2}$.

Let $d$ be a divisor of $2^7 3^{\alpha+1} m$. We
can write $d = 2^{r_1}3^{r_2}t$ where $0 \leq r_1 \leq 7$, $0 \leq r_2 \leq \alpha + 1$ and $t|m$. We now consider the
following three cases depending on $r_1$.
Case 1: Let $0\leq r_1\leq 3$, $0\leq r_2\leq \alpha+1$. Then $G_1=G_2=G_3$, $\dfrac{1}{3^{2\alpha}t^2}\leq G_1\leq 1$ and $G_4=G_5=1$. Hence, $L=3^{\alpha+1}m\left(1-G_1\right)+3\cdot2^k\geq3\cdot2^k.$

Case 2: Let $ r_1=4$, $0\leq r_2\leq \alpha+1$. Then $G_3=G_2=4G_1$, $\dfrac{1}{4\cdot3^{2\alpha}t^2}\leq G_1\leq \dfrac{1}{4}$, $G_5=4G_4$ and $G_4=\dfrac{1}{4}$. Hence, $L=3\cdot2^k.$

Case 3: Let $5\leq r_1\leq 7$, $0\leq r_2\leq \alpha+1$. Then $G_3=4G_2=16G_1$, $\dfrac{1}{16\cdot3^{2\alpha}t^2}\leq G_1\leq \dfrac{1}{16}$, $G_5=4G_4$ and $G_4=\dfrac{1}{16}$. Hence, $L=12\cdot3^{\alpha}m\left(\dfrac{1}{16}-G_1\right)\geq0.$

This proves that $F_{\alpha,m,k}(z)$ is holomorphic at every cusp $c/d$. The weight of $F_{\alpha,m,k}(z)$ is $\ell=\dfrac{1}{2} \left(3^\alpha m+2^k-1\right)$, which is a positive integer and the associated character is given by
 \begin{align*}
\chi_1(\bullet)=\left(\dfrac{(-1)^\ell  2^{4\cdot 3^\alpha m+3\cdot 2^k-4} 3^{(\alpha+1) \left(3^\alpha m+2^k\right)-1} m^{3^\alpha m+2^k}}{\bullet}\right).
\end{align*}

 Thus, $F_{\alpha,m,k}(z)\in M_\ell\left(\Gamma_0(N),\chi\right)$ where $\ell$, $N$ and $\chi$ are as above. Therefore, by Theorem \ref{Ser}, the Fourier coefficients of $F_{\alpha,m,k}(z)$ are almost divisible by $r=2^k$. Due to \eqref{F amk}, this holds for $\overline{b}_{3^{\alpha}m}(n)$ also. This completes the proof of Theorem \ref{b 2k}.
\section{Proof of Theorem \ref{b 3k}}

We define
\begin{align*}
G_{\alpha,m}(z):=\dfrac{\eta^3\left(2^53^{\alpha+1}mz\right)}{\eta\left(2^53^{\alpha+2}mz\right)}.
\end{align*}

Using the binomial theorem, for any integer $k\geq 1$, we have
\begin{align}
G^{3^k}_{\alpha,m}(z)=\dfrac{\eta^{3^{k+1}}\left(2^53^{\alpha+1}mz\right)}{\eta^{3^k}\left(2^53^{\alpha+2}mz\right)}\equiv 1\pmod{3^{k+1}}\label{G 3k}.
\end{align}

Next, we define 
\begin{align*}
H_{\alpha,m,k}(z)&:=\dfrac{\eta(48z)\eta^{3^{\alpha}m}\left(2^33^{\alpha+1}mz\right)\eta^{3^{\alpha}m}\left(2^53^{\alpha+1}mz\right)}{\eta(24z)\eta(96z)\eta^{3^{\alpha}m}\left(2^43^{\alpha+1}mz\right)}G^{3^k}_{\alpha,m}(z)\\
&=\dfrac{\eta(48z)\eta^{3^{\alpha}m}\left(2^33^{\alpha+1}mz\right)\eta^{3^{k+1}+3^{\alpha}m}\left(2^53^{\alpha+1}mz\right)}{\eta(24z)\eta(96z)\eta^{3^{\alpha}m}\left(2^43^{\alpha+1}mz\right)\eta^{3^k}\left(2^53^{\alpha+2}mz\right)}.
\end{align*}

From \eqref{b 3am} and \eqref{G 3k}, we have
\begin{align}
H_{\alpha,m,k}(z)&\equiv\dfrac{\eta(48z)\eta^{3^{\alpha}m}\left(2^33^{\alpha+1}mz\right)\eta^{3^{\alpha}m}\left(2^53^{\alpha+1}mz\right)}{\eta(24z)\eta(96z)\eta^{3^{\alpha}m}\left(2^43^{\alpha+1}mz\right)}\notag\\
&\equiv q^{3\left(3^{2\alpha}m^2-1\right)}\dfrac{f_{48}f_{2^33^{\alpha+1}m}^{3^{\alpha}m}f_{2^53^{\alpha+1}m}^{3^{\alpha}m}}{f_{24}f_{96}f_{2^43^{\alpha+1}m}^{3^{\alpha}m}}\notag\\
&\equiv\displaystyle\sum_{n=0}^{\infty}\overline{b}_{3^{\alpha}m}(n)q^{24n+3\left(3^{2\alpha}m^2-1\right)}\pmod{3^{k+1}}\label{H amk}.
\end{align}

Next, we will prove that $H_{\alpha,m,k}(z)$ is a modular form. Applying Theorem \ref{level}, we find that the level of $H_{\alpha,m,k}(z)$ is $N=2^53^{\alpha+2}mM$, where $M$ is the smallest positive integer such that
\begin{align*}
2^53^{\alpha+2}mM\bigg(\dfrac{-1}{24}+\dfrac{1}{48}+\dfrac{-1}{96}+\dfrac{3^{\alpha}m}{2^33^{\alpha+1}m}+\dfrac{-3^{\alpha}m}{2^43^{\alpha+1}m}\\+\dfrac{3^{k+1}+3^{\alpha}m}{2^53^{\alpha+1}m}
+\dfrac{-3^{k}}{2^53^{\alpha+2}m}\bigg)\equiv 0\pmod{24},
\end{align*}
which gives 
\begin{align*}
8\cdot3^kM\equiv 0\pmod{24}.
\end{align*}
Thus, $M=1$ and the level of $H_{\alpha,m,k}(z)$ is $N=2^5 3^{\alpha+2}m$.

The cusps of $\Gamma_0\left(2^53^{\alpha+2}m\right)$ are given by fractions $c/d$ where $d|2^53^{\alpha+2}m$ and $\gcd(c, d) = 1$. By Theorem \ref{order}, $H_{\alpha,m,k}(z)$ is holomorphic at a cusp $c/d$ if and only if
\begin{align*}
&- \dfrac{\gcd (d,24)^2}{24}+ \dfrac{\gcd (d,48)^2}{48}-\dfrac{\gcd (d,96)^2}{96}+3^{a} m \dfrac{\gcd \left(d,2^33^{\alpha+1} m\right)^2}{2^33^{\alpha+1} m}\\
&- 3^{\alpha} m \dfrac{\gcd \left(d,2^4 3^{\alpha+1} m\right)^2}{2^4 3^{\alpha+1} m}+\left(3^{\alpha} m+3^{k+1}\right) \dfrac{\gcd \left(d,2^53^{\alpha+1} m\right)^2}{2^53^{\alpha+1} m}-3^k \dfrac{\gcd \left(d,2^53^{\alpha+2} m\right)^2}{2^53^{\alpha+2} m}\geq 0.
\end{align*}
Equivalently, $H_{\alpha,m,k}(z)$ is holomorphic at a cusp $c/d$ if and only if
\begin{align*}
L:=3^{\alpha+1}m(-4G_1+2G_2-G_3+4G_4-2G_5+G_6)+3^k(9G_6-1)\geq 0,
\end{align*}
where $G_1=\dfrac{\gcd (d,24)^2}{\gcd \left(d,2^5 3^{\alpha+2} m\right)^2}$, $G_2=\dfrac{\gcd (d,48)^2}{\gcd \left(d,2^5 3^{\alpha+2} m\right)^2}$, $G_3=\dfrac{\gcd (d,96)^2}{\gcd \left(d,2^5 3^{\alpha+2} m\right)^2}$,

 $G_4=\dfrac{\gcd (d,2^3 3^{\alpha+1} m)^2}{\gcd \left(d,2^5 3^{\alpha+2} m\right)^2}$, $G_5=\dfrac{\gcd (d,2^4 3^{\alpha+1} m)^2}{\gcd \left(d,2^5 3^{\alpha+2} m\right)^2}$ and $G_6=\dfrac{\gcd (d,2^5 3^{\alpha+1} m)^2}{\gcd \left(d,2^5 3^{\alpha+2} m\right)^2}$.

Let $d$ be a divisor of $2^53^{\alpha+2} m$. We
can write $d = 2^{r_1}3^{r_2}t$ where $0 \leq r_1 \leq 5$, $0 \leq r_2 \leq \alpha + 2$ and $t|m$. We now consider the
following six cases depending on $r_1$ and $r_2$.

Case 1: Let $0\leq r_1\leq 3$, $0\leq r_2\leq \alpha+1$. Then $G_1=G_2=G_3$, $\dfrac{1}{3^{2\alpha}t^2}\leq G_1\leq 1$ and $G_4=G_5=G_6=1$. Therefore, $L=3^{\alpha+2}m(1-G_1)+8\cdot3^k\geq8\cdot3^k.$

Case 2: Let $0\leq r_1\leq 3$, $ r_2= \alpha+2$. Then $G_1=G_2=G_3$, $\dfrac{1}{3^{2(\alpha+1)}t^2}\leq G_1\leq \dfrac{1}{3^{2(\alpha+1)}}$ and $G_4=G_5=G_6=\dfrac{1}{9}$. Therefore, $L=3^{\alpha+2}m\left(\dfrac{1}{9}-G_1\right)\geq0.$

Case 3: Let $ r_1=4$, $ 0\leq r_2\leq \alpha+1$. Then $G_3=G_2=4G_1$, $\dfrac{1}{4\cdot3^{2\alpha}t^2}\leq G_1\leq \dfrac{1}{4}$, $G_5=G_6=4G_4$ and $G_4=\dfrac{1}{4}$. Hence, $L=8\cdot3^k.$

Case 4: Let $ r_1=4$, $ r_2= \alpha+2$. Then $G_3=G_2=4G_1$, $\dfrac{1}{4\cdot3^{2(\alpha+1)}t^2}\leq G_1\leq \dfrac{1}{4\cdot3^{2(\alpha+1)}}$, $G_5=G_6=4G_4$ and $G_4=\dfrac{1}{36}$. Hence, $L=0.$

Case 5: Let $ r_1=5$, $ 0\leq r_2\leq \alpha+1$. Then $G_3=4G_2=16G_1$, $\dfrac{1}{16\cdot3^{2\alpha}t^2}\leq G_1\leq \dfrac{1}{16}$, $G_6=4G_5=16G_4$ and $G_4=\dfrac{1}{16}$. Therefore, $L=12\cdot3^{\alpha+1}m\left(\dfrac{1}{16}-G_1\right)+8\cdot3^k\geq8\cdot3^k.$

Case 6: Let $ r_1=5$, $ r_2= \alpha+2$. Then $G_3=4G_2=16G_1$, $\dfrac{1}{16\cdot3^{2(\alpha+1)}t^2}\leq G_1\leq \dfrac{1}{16\cdot3^{2(\alpha+1)}}$, $G_6=4G_5=16G_4$ and $G_4=\dfrac{1}{144}$. Therefore, $L=12\cdot3^{\alpha+1}m\left(\dfrac{1}{144}-G_1\right)\geq0.$

This proves that $H_{\alpha,m,k}(z)$ is holomorphic at every cusp $c/d$. The weight of $H_{\alpha,m,k}(z)$ is $\ell=\dfrac{3^\alpha m-1}{2}+3^k$, which is a positive integer and the associated character is given by
 \begin{align*}
\chi_2(\bullet)=\left(\dfrac{(-1)^\ell  2^{4\cdot3^\alpha m+10\cdot 3^k-4} 3^{2 \alpha 3^k+3^\alpha \alpha m+3^\alpha m+3^k-1} m^{3^\alpha m+2\cdot 3^k}}{\bullet}\right).
\end{align*}

 Thus, $H_{\alpha,m,k}(z)\in M_\ell\left(\Gamma_0(N),\chi\right)$ where $\ell$, $N$ and $\chi$ are as above. Therefore, by Theorem \ref{Ser}, the Fourier coefficients of $H_{\alpha,m,k}(z)$ are almost divisible by $r=3^k$. Due to \eqref{H amk}, the same holds for $\overline{b}_{3^{\alpha}m}(n)$ also. This completes the proof of Theorem \ref{b 3k}.

\section{Proof of Theorem \ref{b pk}}

Consider $t=p_1^{a_1}p_2^{a_2}\cdots p_m^{a_m}$, where $p_i$'s are primes. Then we have
\begin{align}
\displaystyle\sum_{n=0}^{\infty}\overline{b}_{t}(n)q^n=\dfrac{f_2f_{t}^{t}f_{4t}^{t}}{f_1f_4f_{2t}^{t}}\label{b t1}.
\end{align}

For a positive integer $i$, we define
\begin{align*}
A_i(z):=\dfrac{\eta^{p_i^{a_i}}(24z)}{\eta\left(24p_i^{a_i}z\right)}
\end{align*}
In view of the binomial theorem, for any integer $k\geq 1$, we have
\begin{align}
A_i^{p_i^k}(z)=\dfrac{\eta^{p_i^{a_i+k}}(24z)}{\eta^{p_i^k}\left(24p_i^{a_i}z\right)}\equiv1\pmod{p_i^{k+1}}\label{Ai pk}.
\end{align}

Define 
\begin{align*}
B_{i,k,t}(z):&=\dfrac{\eta\left(48z\right)\eta^t\left(24tz\right)\eta^t\left(96tz\right)}{\eta\left(24z\right)\eta\left(96z\right)\eta^t\left(48tz\right)}A_i^{p_i^k}(z)\\
&=\dfrac{\eta^{p_i^{a_i+k}-1}\left(24z\right)\eta\left(48z\right)\eta^t\left(24tz\right)\eta^t\left(96tz\right)}{\eta\left(96z\right)\eta^t\left(48tz\right)\eta^{p_i^k}\left(24p_i^{a_i}\right)}.
\end{align*}

From \eqref{b t1} and \eqref{Ai pk}, we arrive at
\begin{align}
B_{i,k,t}(z)&\equiv \dfrac{\eta\left(48z\right)\eta^t\left(24tz\right)\eta^t\left(96tz\right)}{\eta\left(24z\right)\eta\left(96z\right)\eta^t\left(48tz\right)}\notag\\
&\equiv q^{3(t^2-1)}\dfrac{f_{48}f_{24t}^{t}f_{96t}^{t}}{f_{24}f_{96}f_{48t}^{t}}\notag\\
&\equiv \displaystyle\sum_{n=0}^{\infty}\overline{b}_{t}(n)q^{24n+3(t^2-1)}\pmod{p_i^{k+1}}\label{B ikt}.
\end{align}

Next, we show that $B_{i,k,t}(z)$ is a modular form. Applying Theorem \ref{level}, we first estimate
the level of eta quotient $B_{i,k,t}(z)$ . The level of $B_{i,k,t}(z)$ is $N = 96tM$, where
$M$ is the smallest positive integer which satisfies
\begin{align*}
96tM\left(\dfrac{p_i^{a_i+k}-1}{24}+\dfrac{1}{48}+\dfrac{-1}{96}+\dfrac{t}{24t}+\dfrac{-t}{48t}+\dfrac{t}{96t}+\dfrac{-p_i^k}{24p_i^{a_i}}\right)\equiv 0\pmod{24},
\end{align*}
which gives
\begin{align*}
4tM~p_i^k\left(p_i^{a_i}-\dfrac{1}{p_i^{a_i}}\right)\equiv 0\pmod{24}.
\end{align*}
Hence, $M=6$ and $N=2^63^2t.$

The cusps of $\Gamma_0\left(2^63^{2}t\right)$ are given by fractions $c/d$ where $d|2^63^{2}t$ and $\gcd(c, d) = 1$. By Theorem \ref{order}, $B_{i,k,t}(z)$ is holomorphic at a cusp $c/d$ if and only if
\begin{align*}
\left(p_i^{a_i+k}-1\right)\dfrac{\gcd (d,24)^2}{24} +\dfrac{\gcd (d,48)^2}{48}-\dfrac{\gcd (d,96)^2}{96}-p_i^{k} \dfrac{\gcd \left(d,24 p_i^{a_i}\right)^2}{24p_i^{a_i}}\\
+t~\dfrac{ \gcd (d,24 t)^2}{24t}-t~\dfrac{ \gcd (d,48 t)^2}{48t}+\dfrac{\gcd (d,96 t)^2}{96t}\geq 0.
\end{align*}
Equivalently, $B_{i,k,t}(z)$ is holomorphic at a cusp $c/d$ if and only if 
\begin{align*}
L:=-4G_1+2G_2-G_3+4G_4-2G_5+1+4\left(p_i^{k+a_i}G_1-p_i^{k-a_i}G_6\right)\geq0,
\end{align*}
where $G_1=\dfrac{\gcd (d,24)^2}{\gcd (d,96t)^2}$, $G_2=\dfrac{\gcd (d,48)^2}{\gcd (d,96t)^2}$, $G_3=\dfrac{\gcd (d,96)^2}{\gcd (d,96t)^2}$, $G_4=\dfrac{\gcd (d,24t)^2}{\gcd (d,96t)^2}$,

 $G_5=\dfrac{\gcd (d,48t)^2}{\gcd (d,96t)^2}$ and $G_6=\dfrac{\gcd \left(d,24p_i^{a_i}\right)^2}{\gcd (d,96t)^2}$.

Let $d$ be a divisor of $2^63^{2} t$. We
can write $d = 2^{r_1}3^{r_2}p_i^s u$ where $0 \leq r_1 \leq 6$, $0 \leq r_2 \leq  2$, $0\leq s\leq a_i$ and $u|t$ but $p_i\nmid u$. We now consider the
following three cases depending on $r_1$.

Case 1: Let $0\leq r_1\leq 3$, $0\leq r_2\leq 2$. Then $G_1=G_2=G_3=\dfrac{1}{u^2p_i^{2s}}$, $G_4=G_5=1$ and $G_6=\dfrac{1}{u^2}$. Therefore,
\begin{align*}
L=3\left(1-\dfrac{1}{u^2p_i^{2s}}\right)+4\left(\dfrac{p_i^{k+a_i}}{u^2p_i^{2s}}-\dfrac{p_i^{k-a_i}}{u^2}\right)=3\left(1-\dfrac{1}{u^2p_i^{2s}}\right)+4\dfrac{p_i^k}{u^2}\left(\dfrac{p_i^{2a_i}-p_i^{2s}}{p_i^{2s+a_i}}\right).
\end{align*}
Since $s\leq a_i$, we have $L\geq0.$

Case 2: Let $r_1=4$, $0\leq r_2\leq 2$. Then $G_3=G_2=4G_1$, $G_1=\dfrac{1}{4u^2p_i^{2s}}$, $G_5=4G_4$, $G_4=\dfrac{1}{4}$ and $G_6=\dfrac{1}{4u^2}$. Therefore, 
\begin{align*}
L=4\left(\dfrac{p_i^{k+a_i}}{4u^2p_i^{2s}}-\dfrac{p_i^{k-a_i}}{4u^2}\right)=\dfrac{p_i^k}{u^2}\left(\dfrac{p_i^{2a_i}-p_i^{2s}}{p_i^{2s+a_i}}\right)\geq0.
\end{align*}

Case 3: Let $5\leq r_1\leq 6$, $0\leq r_2\leq 2$. Then $G_3=4G_2=16G_1$, $G_1=\dfrac{1}{16u^2p_i^{2s}}$, $G_5=4G_4$, $G_4=\dfrac{1}{16}$ and $G_6=\dfrac{1}{16u^2}$. Hence, 
\begin{align*}
L=\dfrac{3}{4}\left(1-\dfrac{1}{u^2p_i^{2s}}\right)+4\left(\dfrac{p_i^{k+a_i}}{16u^2p_i^{2s}}-\dfrac{p_i^{k-a_i}}{16u^2}\right)=\dfrac{3}{4}\left(1-\dfrac{1}{u^2p_i^{2s}}\right)+\dfrac{p_i^k}{4u^2}\left(\dfrac{p_i^{2a_i}-p_i^{2s}}{p_i^{2s+a_i}}\right)\geq0.
\end{align*}

Therefore, $B_{i,k,t}(z)$ is holomorphic at every cusp $c/d$. The weight of $B_{i,k,t}(z)$ is $\ell=\dfrac{1}{2} \left(p_i^k\left(p_i^{a_i}-1\right)+t-1\right)$, which is a positive integer and the associated character is given by
 \begin{align*}
\chi_3(\bullet)=\left(\dfrac{(-1)^\ell  2^{3 p_i^{a_i+k}-3 p_i^k+4 t-4} 3^{p_i^{a_i+k}-p^k+t-1}t^t \left(p_i^{a_i}\right)^{-p_i^k}}{\bullet}\right).
\end{align*}

 Hence, $B_{i,k,t}(z)\in M_\ell\left(\Gamma_0(N),\chi\right)$ where $\ell$, $N$ and $\chi$ are as above. Therefore, by Theorem \ref{Ser}, the Fourier coefficients of $B_{i,k,t}(z)$ are almost divisible by $r=p_i^k$. Due to \eqref{B ikt}, this holds for $\overline{b}_{t}(n)$ also. Thus, we complete the proof of Theorem \ref{b pk}.
 
 \section{Proof of Theorem \ref{cong2}}\label{pr cong2}
 We recall the following theorem of Ono and Taguchi \cite{ot} on the nilpotency of Hecke operators. 
\begin{theorem}\cite[Theorem 1.3 (3)]{ot}\label{thmot}
Let $n$ be a nonnegative integer and $k$ be a positive integer. Let $\chi$ be a quadratic Dirichlet character of conductor $9\cdot2^a$. Then there is an integer $c \geq 0$ such that for every $f(z) \in M_k(\Gamma_0(9\cdot2^a),\chi) \cap \mathbb{Z}[[q]]$ and every $t\geq 1$,
\begin{align*}
f(z)|T_{p_1}|T_{p_2}|\cdots|T_{p_{c+t}}\equiv 0 \pmod{2^t}
\end{align*}
whenever the primes $p_1, \ldots, p_{c+t}$ are coprime to 6.
\end{theorem}

Now, we apply the above theorem to the modular form $F_{1,1,k}(z)$ to prove Theorem \ref{cong2}.

Taking $\alpha=1$ and $m=1$ in \eqref{F amk}, we have
\begin{align*}
F_{1,1,k}(z)\equiv\displaystyle\sum_{n=0}^{\infty}\overline{b}_{3}(n)q^{24n+24}\pmod{2^{k+1}},
\end{align*}
which yields
\begin{align}
F_{1,1,k}(z):=\displaystyle\sum_{n=0}^{\infty}\mathcal{F}_k\left(n\right)q^{n}\equiv\displaystyle\sum_{n=0}^{\infty}\overline{b}_{3}\left(\dfrac{n-24}{24}\right)q^{n}\pmod{2^{k+1}}\label{F 11k}.
\end{align}

Now, $F_{1,1,k}(z)\in M_{2^{k-1}+1}\left(\Gamma_0(9\cdot 2^7),\chi_4\right)$ for $k\geq 1$ where $\chi_4$ is the associated character (which is $\chi_1$ evaluated at $\alpha=1$ and $m=1$). In view of Theorem \ref{thmot}, we find that there is an integer $u\geq 0$ such that for any $v\geq 1$,
\begin{align*}
F_{1,1,k}(z)\mid T_{q_1}\mid T_{q_2}\mid\cdots\mid T_{q_{u+v}}\equiv 0\pmod{2^v}
\end{align*}
whenever $q_1,\ldots, q_{c+d}$ are coprime to 6. From the definition of Hecke operators, we have that if $q_1,\ldots, q_{u+v}$ are distinct primes and if $n$ is coprime to $q_1\cdots q_{u+v}$,
then
\begin{align}
\mathcal{F}_k\left(q_1\cdots q_{u+v}\cdot n\right)\equiv 0\pmod{2^v}.\label{Fa22}
\end{align}

Combining \eqref{F 11k} and \eqref{Fa22}, we complete the proof of the theorem.

\section{Concluding Remarks}\label{cr}
\begin{enumerate}
\item[(1)] Theorems \ref{b 2k}--\ref{b pk} give the arithmetic densities of $\overline{b}_t(n)$ for odd $t$. But it is not possible to study the arithmetic density of $\overline{b}_t(n)$ for even $t$ using the similar techniques. It would be interesting to study the arithmetic density of $\overline{b}_t(n)$ when $t$ is even.

\item[(2)]Some congruences for $\overline{b}_5(n)$ modulo powers of 5 have been proved in \cite{bb}. We encourage the readers to find new Ramanujan-type congruences for $\overline{b}_t(n)$ modulo powers of primes for different values of $t$. 

\item[(3)] Bandyopadhyay and Baruah \cite{bb} deduced several arithmetic identities for $\overline{b}_5(n)$. They also found relations connecting $\overline{a}_t(n)$, $\overline{b}_t(n)$ and $c_t(n)$. It will be desirable to find similar arithmetic identities involving $\overline{b}_t(n)$ for other values of $t$ also. 
\end{enumerate}

\section{Acknowledgement}
The author is extremely grateful to his Ph.D. supervisor, Prof. Nayandeep Deka Baruah, for his guidance and encouragement. The author is indebted to Prof. Rupam Barman for many helpful comments and suggestions. The author was partially supported by the Council of Scientific \& Industrial Research (CSIR), Government of India under the CSIR-JRF scheme (Grant No. 09/0796(12991)/2021-EMR-I). The author thanks the funding agency.

\end{document}